\documentclass[12pt]{amsart}

\usepackage[parfill]{parskip}
\usepackage[toc,page]{appendix}
\usepackage[normalem]{ulem}
\usepackage{psfrag}
\usepackage{graphicx}
\usepackage{xcolor}
\usepackage{pinlabel}
\usepackage{enumitem}
\usepackage{amsmath}
\usepackage{amssymb}
\usepackage{amscd}
\usepackage{picinpar}
\usepackage[pagewise]{lineno}

\usepackage{pb-diagram}
\usepackage{graphicx}
\usepackage{wrapfig}
\usepackage{xspace}
\usepackage{hyperref}
\usepackage{xcolor}
\usepackage{url}
\usepackage{float}

\usepackage{cite}

\usepackage{comment}

\def\R{\mathbb{R}}

\def\T{\mathcal{T}}

\newtheorem{theorem}{Theorem}[section]
\newtheorem{lemma}[theorem]{Lemma}
\newtheorem{proposition}[theorem]{Proposition}

\newtheorem{definition}[theorem]{Definition}

\newtheorem*{theorem*}{Theorem}

\newtheorem{prop}[theorem]{Proposition}

\newtheorem{rmk}[theorem]{Remark}

\def\k{\mathcal{K}}

\def\ulu{\underline{u}}

\title{Approaching the prescribed Gaussian curvature by discrete conformality}
\author{Ziran Liu}
\address{Shanghai Institute for Mathematics and Interdisciplinary Sciences}
\email{zliu@simis.cn}

\author{Tianqi Wu}
\address{Clark University}
\email{mike890505@gmail.com}
\date{}  
\begin{document}



\begin{abstract}
We propose a discrete approach for approximating solutions to the prescribed Gaussian curvature problem in two-dimensional manifolds, based on the notion of discrete conformality. Our approach provides an efficient numerical method to compute the solution by minimizing a convex functional. 
\end{abstract}

\maketitle

\section{Introduction}\label{intro}
The problem of the prescribed Gaussian curvature problem is classical and extensively studied in differential geometry. Suppose $(M, g)$ is a $2$-dimensional Riemannian manifold. Given a smooth function $\tilde {\kappa}: M \rightarrow \mathbb{R}$, the prescribed Gaussian curvature problem asks whether there exists a Riemannian metric $\tilde{g}$ on $M$ such that the Gaussian curvature of $(M,\tilde{g})$ is exactly $\tilde{\kappa}$. By letting $\tilde{g} = e^{2\tilde u}g$ for some smooth function $\tilde u: M \rightarrow \R$, it suffices to solve the equation
\begin{equation} \label{conformal factor u}
    \Delta_{g} \tilde u + \kappa = -e^{2\tilde u}\tilde{\kappa},
\end{equation}
where $\kappa$ is the Gaussian curvature associated with the metric $g$, and $\Delta_{g}$ is the Laplace-Beltrami operator with respect to $g$. The prescribed Gaussian curvature problem has been extensively investigated, see \cite{berger1971riemannian}, \cite{kazdan1974curvature}, \cite{kazdan1975scalar}, 
for examples. 

This paper aims to give a discrete approach to the prescribed Gaussian curvature problem by approximating the solution $u$ to the equation (\ref{conformal factor u}), for a class of $2$-dimensional manifolds. A key ingredient is the theory of \textit{discrete conformality}, which was developed in Luo \cite{luo2004combinatorial} and Bobenko-Pinkall-Springborn \cite{bobenko2015discrete} for \textit{triangle meshes} based on \textit{vetex scalings}.

Recall a well-known result for the existence of the solution to prescribed Gaussian curvature problem on surfaces with $genus > 1$ (Theorem 1.1 in \cite{berger1971riemannian}, also see \cite{kazdan1974curvature}).

\begin{prop} \label{sufficient condition for g > 1}
    Given a smooth function $\tilde{\kappa}(x)<0$ on $(M,g)$, there exists a unique function $\tilde u$ on $M$ such that $e^{2\tilde u}g$ has the Gaussian curvature $\tilde{\kappa}$. Furthermore, such $\tilde  u$ is smooth on $M$.
\end{prop}

This paper will describe a discrete approach to compute the conformal factor $u$ in Proposition \ref{sufficient condition for g > 1}.

\subsection{Discrete conformality for polyhedral surfaces}

Let $S$ be a closed topological surface and $\T$ be a triangulation of $S$. Let $F=F(\T)$ denote the set of faces ($2$-simplexes) of $\T$, $E = E(\T)$ denote the set of edges ($1$-simplexes) of $\T$, and $V = V(\T)$ denote the set of vertices of $\T$. A \textit{polyhedral metric} (PL metric) on $\T$ could be identified as an edge length function $l\in\mathbb R^E_{>0}$. We call the pair $(\T,l)$ a \textit{polyhedral surface}. 
Given $(T,l)$ and a triangle $\triangle ijk$ with the vertices $i, j, k \in V$, let $\theta^{i}_{jk}$ denote the inner angle at vertex $i$ in the Euclidean triangle $\triangle ijk$. 

\begin{definition}\label{d curvature}
The \textit{discrete curvature} is a function
$$
K: V \rightarrow (-\infty, 2\pi),
$$ 
such that 
$$ 
K_i =2\pi- \sum_{jk\in E:\triangle ijk\in F}\theta^i_{jk}.
$$
\end{definition}

\begin{definition}[Discrete Conformality by Vertex Scaling, \cite{luo2004combinatorial}]

\label{discrete conformality with fixed triangulation}
Given a triangulation $\T$ on $(S,V)$, two Euclidean polyhedral surfaces $(\T,l)$ and $(\T,l')$ are {\emph{discretely conformally equivalent}}, if for some $u \in \R^{V(\T)}$, 
\begin{equation}
\label{Euclidean discrete conformal}
l'_{ij}=e^{\frac{1}{2}(u_i+u_j)}l_{ij}
\end{equation}
for any edge $ij\in E(\T)$. 
We call such a $u$ by \emph{discrete conformal factor}, and denote $l'=u*l$ if equation (\ref{Euclidean discrete conformal}) holds.
\end{definition}

\subsection{Triangulations with mixed background constant curvatures}

Suppose $\k \in\mathbb R^{F(\T)}_{<0}$ and $l\in\mathbb R^{E(\T)}_{>0}$. Whenever well-defined, for any $\sigma\in F(\T)$ we denote  

(a) $g(\sigma,\k,l)$ as a Riemmanian metric on $\sigma$ such that $(\sigma,g(\sigma, \k, l))$ is a geodesic triangle of the constant Gaussian curvature $\k(\sigma)$ such that the lengths of an edge $e$ is 
$$
\frac{2}{-\k(\sigma)}\sinh^{-1}\left(\frac{-\k (\sigma)}{2}l(e)\right),
$$

(b)
$\theta^i_{jk}=\theta^i_{jk}(\k,l)$ as the value of inner angle at the vertex $i$ in the geodesic triangle 
$$
(\triangle ijk,g(\triangle ijk,\k,l)),
$$
and

(c)
$K=K(\k,l)\in\mathbb R^V$ as a generalized discrete curvature such that 
$$ 
K_i =2\pi- \sum_{jk\in E:\triangle ijk\in F}\theta^i_{jk}.
$$

\begin{rmk}
    Note that in item (c), the angle $\theta^i_{jk}$'s at the same vertex $i$ may not share the same metric as in the definition \ref{d curvature} of the \emph{discrete curvature}. 
\end{rmk}

As observed in Bobenko-Pinkall-Springborn \cite{bobenko2015discrete}, if $\k(\sigma)=-1$, then $(\sigma,g(\sigma,\k,u*l))$ is discretely conformal to $(\sigma,g(\sigma,\k,l))$ as hyperbolic geodesic triangles with the discrete conformal factor $u$, in the notion of hyperbolic discrete conformality introduced in \cite{bobenko2015discrete}. 
We can naturally extend their notion of hyperbolic discrete conformality to geodesic triangles of constant negative curvature $\k(\sigma)$. 
For fixed $\k\text{ and }l$, we denote $K(\k,u*l)$ by $K(u)$ for convenience.

Then by the variational principle developed by Luo \cite{luo2004combinatorial} and Bobenko-Pinkall-Springborn \cite{bobenko2015discrete}
$$
\mathcal F(u)=\int \sum_{i\in V} K_i(u)du_i
$$
is well-defined and locally strictly convex on the domain where $K(u)$ is well-defined. 

Furthermore, 
 Bobenko-Pinkall-Sprinborn \cite{bobenko2015discrete} gave a simple and explicit formula extending $\mathcal F$ to a globally convex functional $\tilde{\mathcal F}$ on $\mathbb R^V$. So $K(u)=0$ has at most one solution, and 
  can be computed efficiently by minimizing a globally convex functional $\tilde {\mathcal F}$ on $\mathbb R^V$.

One can regard $\mathcal{K} \in\mathbb R^{F(\T)}_{<0}$ as a discrete approximation of some smooth Gaussian curvature function $\tilde \kappa \in \R^S$. Thus, when $K(u) = 0$, one may regard $\cup_\sigma(\sigma,g(\sigma,\k,u*l))$ as an approximation of a Riemannian surface with the prescribed Gaussian curvature $\tilde \kappa$. Meanwhile, $u \in \R^{V(\T)}$ may be regarded as a discrete approximation of the smooth conformal factor $u \in \R^S$ satisfying that $e^{2u}g$ has the Gaussian curvature $\tilde \kappa$.

\subsection{Geodesic triangulation, and discrete approximation of the prescribed curvature problem}
For triangulation $\T$, if any edge in $\T$ is a geodesic arc on $S$, we call such $\T$ a \emph{geodesic triangulation}.

Suppose that $\T$ is a geodesic triangulation of $(S,g)$ and $V, E, F$ are defined the same way as in section \ref{intro}. The edge length function $l\in\mathbb R^{E(\T)}_{>0}$ is measured on $(S,g)$, and the polyhedral surface $(\T, l)$ is regarded as an approximation of $(S,g)$. 

Suppose $\tilde \kappa(x)<0$ is a smooth function on $S$, by Theorem 1.1, there exists a smooth function $ \tilde u$
such that $e^{2\tilde u}g$ has the Gaussian curvature $\tilde\kappa$. 

We will find a discrete approximation of $\tilde u$ using the information of $(\T, l)$ and $\tilde\kappa$.
For all $\sigma\in F(\T)$, we let $\k(\sigma)=\tilde \kappa(x)$ for some $x\in \sigma$. Then by minimizing a convex functional on $\mathbb R^V$, we expect to find a unique discrete conformal factor $u\in\mathbb R^V$ such that $K(\k,u*l)=0$ (i.e. $K(u) = 0$). Such a discrete conformal factor $u$ should be an approximation of $\tilde u$ in the sense of item $(b)$ in Theorem \ref{high genus main thm}.

We need some regularity of the triangulation to prove the convergence. Recall from section \ref{intro}, for $i, j, k \in V$, $jk \in E$ and $ijk \in F$, we let $\theta^{i}_{jk}$ be the value of the inner angle at vertex $i$ in the triangle $\triangle ijk$. 

\begin{definition}\label{epsilon regular mesh}
A polyhedral surface $(\T, l)$ is called \emph{$\epsilon$-acute} if
  $$
  \theta^i_{jk}\leq\frac{\pi}{2}-\epsilon
  $$ 
  for all inner angle.
\end{definition}

In this paper, if $x\in\mathbb R^A$ is a vector for some finite set $A$, we use
$|x|$ to denote the infinity norm of $x$, i.e., $|x|=|x|_\infty=\max_{i\in A}|x_i|$.
Now we state our main theorem.
\begin{theorem}\label{high genus main thm}
Suppose $(S,g)$ is a connected closed orientable smooth Riemannian surface with $genus > 1$, $\tilde \kappa(x)<0$ is a smooth function on $S$, and $\tilde u$ is the smooth function on $S$ such that $e^{2\tilde u}g$ has Gaussian curvature $\tilde\kappa(x)$.
Assume $\T$ is a geodesic triangulation of $(S,g)$, and $l\in\mathbb R^{E(T)}_{>0}$ denotes the edge length in $(S,g)$. For any triangle $\sigma$ in $F(\T)$, let $\k(\sigma)=\tilde \kappa(x)$ for some
$x\in\sigma$.
Then for any $\epsilon>0$, there exist constants $\delta=\delta(S,g,\tilde \kappa,\epsilon)>0$ and $C=C(M,g,\tilde \kappa,\epsilon)$ such that if $(\T, l)$ is $\epsilon$-acute and $|l|<\delta$, then
\begin{enumerate}[label=(\alph*)]
\item there exists a unique discrete conformal factor $u\in \mathbb R^{V(\T)}$, such that $K(\mathcal K,u*l)=0$, and
\item$\big|u-\tilde u|_{V(\T)}\big|\leq C|l|.$
\end{enumerate}
\end{theorem}
The main idea of the proof of Theorem 1.4 follows \cite{wu2020convergence}. We describe a key discrete elliptic estimate in Section 2 and the formula of $\partial K/\partial u$ in Section 3. The proof of Theorem 1.4 is given in Section 4.

\section{Discrete calculus on graphs}
\label{section 2}
Assume $G=(V,E)$ is an undirected connected simple graph. Let $ij$ denote an edge in $E$ with endpoints $i$ and $j$, and $i \sim j $ represent that the vertices $i$ and $j$ are connected by an edge $ij\in E$. Let $\mathbb R^E$ and $\mathbb R^E_A$ be vector spaces of dimension $|E|$ such that
\begin{enumerate}[label=(\alph*)]
	\item a vector $x\in \mathbb R^E$ is represented symmetrically, i.e., $x_{ij}=x_{ji}$, and
	\item a vector $x\in\mathbb R^E_A$ is represented anti-symmetrically, i.e., $x_{ij}=-x_{ji}$.
\end{enumerate}
A vector $x\in \mathbb R^E_A$ is called a \emph{flow} on $G$.
An \emph{edge weight} $\eta$ on $G$ is a vector in $\mathbb R^E$. Given an edge weight $\eta$, the \emph{gradient} $\nabla f=\nabla_\eta f$ of a vertex function $f\in\mathbb R^V$ is a flow in $\mathbb R^E_A$ defined as
$$
(\nabla f)_{ij}=\eta_{ij}(f_j-f_i)
$$
where $f_i = f(i).$
Given a flow $x\in\mathbb R^E_A$, its \emph{divergence} $div(x)$ is a vector in $\mathbb R^V$ such that
$$
div(x)_i=\sum_{j\sim i}x_{ij}.
$$
Given an edge weight $\eta$, the associated \emph{Laplacian} $\Delta=\Delta_\eta:\mathbb R^V\rightarrow\mathbb R^V$ is defined as
$\Delta f=div(\nabla_\eta f)$, i.e.,
$$
(\Delta f)_i=\sum_{j\sim i}\eta_{ij}(f_j-f_i).
$$
A Laplacian $\Delta_{\eta}$ is a linear transformation on $\mathbb R^V$, and can be identified as a $|V|\times |V|$ symmetric negative semi-definite matrix.

Now we introduce the notion of $C$-isoperimetry for a graph $G=(V,E)$ associated with a positive vector $l\in\mathbb R^E_{>0}$.
Given any $V_0\subset V$, denote
$$
\partial V_0=\{ij\in E:i\in V_0,j\notin V_0\},
$$
and then define the $l$-\emph{perimeter} of $V_0$ and the $l$-\emph{area} of $V_0$ as
\begin{equation*}
	|\partial V_0|_l=\sum_{ij\in \partial V_0}l_{ij}\quad\text{ and }\quad
	|V_0|_l=\sum_{i,j\in V_0,ij\in E}l_{ij}^2
\end{equation*}
respectively.

For a constant $C>0$,
such a pair $(G,l)$ is called $C$-\emph{isoperimetric} if for any $V_0\subset V$
\begin{equation*}
\min\{|V_0|_l,|V|_l-|V_0|_l\}\leq C\cdot|\partial V_0|_l^2.
\end{equation*}
We will see, from part (b) of Lemma \ref{isoperimetric}, that a uniform $C$-isoperimetric condition is satisfied by regular geodesic triangulations approximating a closed smooth surface.
The following discrete elliptic estimate is a key ingredient of the proofs of our main theorems. 
\begin{lemma}[Wu-Zhu \cite{wu2020convergence}]
	\label{estimate for divergence operator}
Assume $(G,l)$ is $C_1$-isoperimetric, and $x\in\mathbb R^E_A,\eta\in\mathbb R^E_{>0},C_2>0,C_3>0$ are such that
\begin{enumerate}[label=(\roman*)]
\item $|x_{ij}|\leq C_2 l_{ij}^2$ for any $ij\in E$, and

\item $\eta_{ij}\geq C_3$ for any  $ij\in E$.
\end{enumerate}
Then
\begin{equation*}
|\Delta^{-1}_\eta\circ div (x))|\leq \frac{4C_2\sqrt{C_1+1}}{C_3}|l|\cdot|V|_l^{1/2}.
\end{equation*}
Furthermore, if $y\in\mathbb R^V$ and $C_4>0$ and $D\in\mathbb R^{V\times V}$ is a non-zero diagonal matrix such that
$$
|y_i|\leq C_4 D_{ii}|l|\cdot|V|_l^{1/2}
$$
for any $i\in V$, then
\begin{equation*}
|(D-\Delta_\eta)^{-1} (div(x)+y)|\leq
\left(C_4+\frac{8C_2\sqrt{C_1+1}}{C_3}\right)|l|\cdot|V|_l^{1/2}.
\end{equation*}
\end{lemma}

\section{Differential of the discrete curvatures and angles}
\label{differential of discrete curvatures and angles}
The following proposition can be easily derived from Proposition 6.1.7 in \cite{bobenko2015discrete} or
Lemma 3.4 in \cite{wu2020convergence}.
\begin{proposition}
	\label{Discrete curvature equation for hyperbolic case}
Given $(\T,l)$, $\k\in\mathbb R^{F(\T)}_{<0}$, and $u\in\mathbb R^{V(\T)}$  such that $K(u)$ is well-defined,
denote
$$
\hat\theta^i_{jk}(u)=\frac{1}{2}(\pi+\theta^i_{jk}(u)-\theta^j_{ik}(u)-\theta^k_{ij}(u))
$$
and
$$
\lambda_{ij,k}=
\frac{\k(\triangle ijk)^2\cdot(u*l)_{ij}^2}{\k(\triangle ijk)^2\cdot(u*l)_{ij}^2+4}
$$
for all triangle $\triangle ijk\in F(\T)$.
Then
\begin{equation*}
\frac{\partial K}{\partial u}(u)=D(u)-\Delta_{\eta(u)}
\end{equation*}
where
$$
\eta_{ij}(u)=
\frac{1}{2}\cot\hat\theta^k_{ij}(u)(1-\lambda_{ij,k})+
\frac{1}{2}\cot\hat\theta^{k'}_{ij}(u)(1-\lambda_{ij,k'}),
$$
and
$D=D(u)$ is a diagonal matrix such that
$$
D_{ii}(u)=\sum_{\triangle ijk\in F}
(
\cot\hat\theta^k_{ij}(u)\lambda_{ij,k}+
\cot\hat\theta^{k'}_{ij}(u)\lambda_{ij,k'}).
$$

\end{proposition}

\section{Proof of the Main Theorems}
\label{section 4}
In this section we first introduce two geometric lemmas, and then prove Theorem \ref{high genus main thm} in Subsections \ref{the case of high genus}.
\subsection{Two geometric lemmas}
The following lemma was first proved in \cite{gu2019convergence}, and was also proved in \cite{wu2020convergence}.
\begin{lemma}
	\label{cubic estimate}
Suppose $(S,g)$ is a closed Riemannian surface, and $ u\in C^{\infty}(S)$. Then there exists $C=C(S,g,u)>0$ such that for any $x,y\in S$,
$$
|d_{e^{2 u}g}(x,y)-e^{\frac{1}{2}(u(x)+ u(y))}d_g(x,y)|\leq Cd_g(x,y)^3,
$$
where $d_{e^{2 u}g}(x,y)$ and $d_g(x,y)$ are the distance metrics from $x$ to $y$ under ${e^{2 u}g}$ and $g$ respectively.
\end{lemma}

The following Lemma \ref{isoperimetric} is similar to Lemma 4.4 in \cite{wu2020convergence} and can be proved almost in the same way.
\begin{lemma}
	\label{isoperimetric}
Suppose $(S,g)$ is a closed Riemannian surface,
    and $\T$ is a geodesic triangulation of $(S,g)$ with $V, E$ the set of vertices and edges. Let $l\in\mathbb R^{E(\T)}$ be the geodesic lengths of the edges and
assume that $(\T, l)$ is $\epsilon$-acute.
\begin{enumerate}[label=(\alph*)]

\item Given $ \tilde u\in C^{\infty}(S)$, there exists a constant $\delta=\delta(S,g,\tilde u,\epsilon)>0$ such that if $|l|<\delta$ then
there exists a geodesic triangulation $\T'$ in $(S,e^{2\tilde u}g)$ such that $V(\T')=V(\T)$, and $\T'$ is homotopic to $\T$ relative to $V(\T)$. Furthermore, $(S, V, \bar l)$ is $\frac{1}{2}\epsilon$-acute where $\bar l\in\mathbb R^{E(\T')}$ denotes the geodesic lengths of the edges of $\T'$ in $(S,e^{2\tilde u}g)$.
\item There exist constants $\delta=\delta(S,g,\epsilon)$ and $C=C(S,g,\epsilon)$ such that if $|l|<\delta$,
then $(\T, l)$ is $C$-isoperimetric.
\end{enumerate}
\end{lemma}

\subsection{Proof of Theorem \ref{high genus main thm}}
\label{the case of high genus}
For a constant $\epsilon>0$, we assume that $(\T, l)$ is $\epsilon$-acute and $|l|\leq \delta$ where $\delta=\delta(S,g,\tilde \kappa,\epsilon)<1$ is a sufficiently small constant to be determined. By Lemma \ref{isoperimetric},
if $\delta$ is sufficiently small there exists a geodesic triangulation $\T'$ of $(S,e^{2\tilde u}g)$ homotopic to $\T$ relative to $V(\T)=V(\T')$.
Let $\bar l\in\mathbb R^{E(\T)}\cong R^{E(\T')}$ denote the geodesic lengths of edges of $\T'$ in $(S,e^{2\tilde u}g)$. 

For simplicity, we will frequently use the notion $A =O(B)$ to denote that if $\delta=\delta(S,g,\tilde \kappa,\epsilon)$ is sufficiently small, then $|A|\leq C\cdot |B|$ for some constant $C=C(S,g,\tilde\kappa,\epsilon)$. For example, we have that
\begin{enumerate}[label=(\alph*)]
\item $l_{ij}=O(l_{jk})$ for any $\triangle ijk\in F(\T)$, 

\item $(\tilde u*l)_{ij}=O(l_{ij})$,

\item $\bar l_{ij}=O(l_{ij})$, and

\item $l_{ij}=O(\frac{2}{-\k(\triangle ijk)}\sinh^{-1}(\frac{-\k(\triangle ijk)}{2}l_{ij}))$.
\end{enumerate}

The remaining of the proof is divided into two steps.
\begin{enumerate}[]
\item Firstly we show that 
$$
(\tilde u*l)_{ij}-\bar l_{ij}=O(l_{ij}^3)
$$
and
$$
K(\tilde u)=div(x)+y
$$
for some $x\in\mathbb R^E_A$ and $y\in\mathbb R^V$ such that $x_{ij}=O(l_{ij}^2)$ and $y_i=O(l_{ij}^3)$.

\item Secondly, we construct a smooth path $\underline{u}(t):[0,1]\rightarrow\mathbb R^V$ with $u(0)=\tilde u|_{V}$ such that the following identity holds
    $$
    K(\ulu(t))=(1-t)K(\tilde u).
    $$ 
    Furthermore, we show that $|\ulu'(t)|=O(|l|)$, and then $K(\ulu(1))=0$ and $\ulu(1)-\tilde u=O(|l|)$.
\end{enumerate}
Therefore the vertex function $\ulu(1)$ is the wanted discrete conformal factor for $(\T, l)$ in Theorem \ref{high genus main thm}.

\subsubsection{Part 1}
\label{hyperbolic 1}
By Lemma \ref{isoperimetric}, $(S, V,\bar{l})$ is $\frac{1}{2}\epsilon$-acute if $\delta$ is sufficiently small. For simplicity we denote $\tilde u|_{V(\T)}$ by $\tilde u$. By Lemma \ref{cubic estimate}, we have that
\begin{equation*}
\label{11}
\bar l_{ij}-(\tilde u*l)_{ij}=O (l_{ij}^3).
\end{equation*}
Notice that 
$$
\tilde \kappa(x)-\k(\triangle ijk)=O(l_{ij})
$$ 
for all $\triangle ijk\in F(\T)$ and $x\in \triangle ijk$. 
Denote $\tilde \theta^i_{jk}$ as the inner angle of the geodesic triangle $(\triangle ijk,e^{2\tilde u}g)$ in $\T'$.
By standard estimates on geodesic triangles with bounded Gaussian curvature, we have that
\begin{equation}
\label{flow in hyperbolic}
\alpha^i_{jk}:=\tilde\theta^i_{jk}-\theta^i_{jk}(\tilde u)=O(l_{ij}^2),
\end{equation}
and
\begin{equation}
\label{residue in hyperbolic}
\alpha^i_{jk}+\alpha^j_{ik}+\alpha^k_{ij}=O(l^3_{ij}).
\end{equation}
So $(\T,\tilde u*l)$ is $\frac{1}{3}\epsilon$-acute if $\delta$ is sufficiently small.
Let $x\in\mathbb R^E_A$ and $y\in\mathbb R^V$ be such that
\begin{equation}
    \label{definition of hyperbolic equation}
x_{ij}=\frac{\alpha^i_{jk}-\alpha^j_{ik}}{3}+\frac{\alpha^i_{jk'}-\alpha^j_{ik'}}{3}\quad
\text{ and }
\quad
y_i=\frac{1}{3}\sum_{jk:\triangle ijk\in F(\T)}(\alpha^i_{jk}+\alpha^j_{ik}+\alpha^k_{ij})
\end{equation}
where $\triangle ijk$ and $\triangle ijk'$ are adjacent triangles. Then we have the following identity
\begin{equation*}
div(x)_i+y_i=K_i(\tilde u).
\end{equation*}
 By (\ref{flow in hyperbolic}), (\ref{residue in hyperbolic}), (\ref{definition of hyperbolic equation}) and the fact that any vertex $i\in V(\T')$ has at most $2\pi/(\epsilon/2)$ neighbors, we have
\begin{equation}
\label{21}
x_{ij}=O(l_{ij}^2),
\end{equation}
and
\begin{equation}
\label{22}
y_i=O(l_{ij}^3).
\end{equation}
\subsubsection{Part 2}
Let
$$
\tilde\Omega=\{u\in\mathbb R^V:u*l\text{ satisfies the triangle inequalities and $(S, V, u*l)$ is $\frac{\epsilon}{4}$-acute}\}
$$
and
$$
\Omega=\{u\in\tilde\Omega:|u-\tilde u|\leq1\text{, $(S, V, u*{l})$ is $\frac{\epsilon}{5}$-acute}\}.
$$

Since $(\T,\tilde u*l)$ is $\frac{1}{3}\epsilon$-acute, $\tilde u$ is in the interior of $\Omega$.
Now consider the following ODE on $int(\tilde\Omega)$,
\begin{equation}
\label{ode for high genus case}
	\left\{
	\begin{array}{lcl}
	u'(t)=(\Delta_{\eta(u)}-D(u))^{-1} K(\tilde u)=(\Delta_{\eta(u)}-D(u))^{-1}(div(x)+y)\\
	u(0)=\tilde u
	\end{array}
	\right.,
\end{equation}
where $D(u)$ and $\eta(u)$ are defined as in Proposition \ref{Discrete curvature equation for hyperbolic case}. For any triangle $\triangle ijk$ and $u\in\tilde\Omega$, we have
\begin{equation}
\label{24}
D_{ii}(u)\geq \epsilon'l_{ij}^2,
\quad
\text{ and }
\quad \eta_{ij}(u)\geq\epsilon'
\end{equation}
for some constant $\epsilon'=\epsilon'(S,g,\tilde \kappa,\epsilon)>0$.

The right-hand side of equation (\ref{ode for high genus case}) is a smooth function of $u$, so the ODE (\ref{ode for high genus case}) has a unique solution $u(t)$ satisfying
$$\frac{dK(\ulu(t))}{dt} = \frac{\partial K}{\partial u}(u)\ulu'(t) 
= 
(D(u)-\Delta_{\eta(u)})
(\Delta_{\eta(u)}-D(u))^{-1}
K(\tilde u)
=-K(\tilde{u}).$$
Therefore
$$
K(\ulu(t))=(1-t)K(\tilde u).
$$
Assume the maximum existence interval of $u(t)$ in $\Omega$ is $[0,T_0)$ where $T_0\in(0,\infty]$.
By Lemma \ref{isoperimetric}, when $\delta$ is sufficiently small, $(\T, l)$ is $C$-isoperimetric for some constant $C=C(S, g, \epsilon)$. Then for any $u\in\Omega$, $(S, V, u*l)$ with triangulation $\T$ is $(e^{4(|\tilde u|+1)}C)$-isoperimetric.
It is not difficult to see
$$
|V|_l=O(1)
$$
and
$$
1=O(|V|_l).
$$
Then by Lemma \ref{estimate for divergence operator} and equation (\ref{21}),(\ref{22}),(\ref{24}), for any $t\in[0,T_0)$
\begin{equation}
|u'(t)| =O(|l| \cdot|V|_l^{1/2})=O(|l|).
\end{equation}
If $T_0\leq1$,
\begin{equation}
|\ulu(T_0)-\tilde u|=O(|l|)\quad\text{ and }\quad\theta^i_{jk}(\ulu(T_0))-\theta^i_{jk}(\tilde u)=O(|l|),
\end{equation}
and then $\ulu(T_0)\in \text{int}(\Omega)$ if $\delta$ is sufficiently small. But this contradicts with the maximality of $T_0$. So $T_0>1$ and
 $K(u(1))=0$ and
$|\ulu(1)-\tilde u|=O(|l|).$

\bibliographystyle{plain}
\bibliography{kc2gkc2}

\end{document}